\documentclass[12pt]{article}
\usepackage{amsmath, amsfonts, amsthm}
\usepackage{geometry}

\usepackage{color}

\newcommand{\bbR}{\mathbb{R}}
\newcommand{\bS}{\overline{S}}

\newif\ifno
\nofalse

\title{{\normalsize\tt\hfill\jobname.tex}\\
On robustness of discrete time optimal filters
}
\author{M.~L.~Kleptsyna\footnote{Universit\'e du Maine, Le Mans, France,
e-mail: Marina.Kleptsyna at univ.lemans.fr} $\;\;$  \& $\;$
A.~Yu.~Veretennikov\footnote{University of Leeds, UK; National Research University Higher School of Economics, \&
Institute of Information Transmission Problems,  Moscow, Russia,
e-mail: A.Veretennikov at leeds.ac.uk} 
}

\begin{document}

\newtheorem{Theorem}{Theorem}
\newtheorem{Lemma}{Lemma}
\newtheorem{Corollary}{Corollary}
\newtheorem{Remark} {Remark}
\newtheorem{Proposition}{Proposition}

\maketitle

\begin{abstract}
A new result on stability of an optimal nonlinear filter for a Markov chain with respect to small perturbations on every step is established. An  exponential recurrence of the signal is assumed.
\end{abstract}

\section{Introduction}
Stability of optimal filters is a topical research area in the last three or even more decades.  In this direction, a lot has been understood and achieved under the ``uniform ergodicity'' assumptions due to the method by Atar and Zeitouni (see \cite{AZ:97}) based on the Birkhoff metric (also known as projective or Hilbert metric). This method under such assumptions guarantees an exponential rate, with which the optimal filter algorithms ``forgets'' wrong -- or unspecified -- initial conditions. The method has been extended to the ``non-uniform ergodic'' case (see \cite{KV:08}) by combining the application of Birkhoff metric with a modified version of the coupling method, which leaded to exponential or polynomial rates, with which the algorithm  may ``forget'' wrong initial data. However, unspecified initial distributions is not the only option for an unspecified model. Small errors on each step of the algorithm (in discrete time case) is one more possibility to perturbe a model.  In the ``uniform ergodic'' case it was also tackled in the literature (see \cite{LGO:04}). The ``non--uniform ergodic'' case is still waiting for its investigation and our goal here is to attack this problem. In our setting only ``uniformly small'' errors are allowed and conditions on the densities of the noise both in the signal and in the observations look rather strict, so that new studies will be required to weaken conditions so as to include a wider class of processes.

The setting described earlier is not only insteresting as such: it may also serve as a base for studying unspecified models with an unknown parameter. In such models, observations should allow to estimate the parameter. Once the estimator is, at least, consistent, there is a hope that the filtering algorithm for a model with an estimate instead of the ``true parameter'' may be close enough to the exact model. Hence, the previous studies could be applied. This programme -- again in the ``uniform'' case -- was realised in \cite{Ana1, Ana2}, where it was {\em assumed} that the estimator satisfies certain large deviation conditions. However, in many examples ``non-uniform''
conditions are more natural. Hence, a large part of the problem remains open and requires further investigations. We restrict ourselves to the case of exponential recurrence and exponential moments for the signal and postpone other cases till further research.

The paper consists of three sections: Introduction, Setting and main result, which also includes some auxiliary results; Proof of auxiliary lemmata; Proof of main result -- the Theorem \ref{Thm1}.

\section{Setting and Main Result}\label{sett}
We consider the following model, with a non-observed (Markov) state process $\{X_{n}, n\ge 0\}$ and an observation process
$\{Y_{n}, n\ge 0\}$, taking value in $\mathbb{R}^{d}$ and $\mathbb{R}^{\ell}$ respectively.
We assume that the state sequence $\{X_{n}, n\ge 0\}$ is defined as a  homogeneous Markov chain with transition probability kernel  $Q(x,dx')$, \textit{i.e.}:
\begin{equation}\label{kernel}
\mathbb{P}[X_{n}\in dx' | X_{0:n-1}]|_{X_{n-1}= x} =\mathbb{P}[X_{n}\in dx' | X_{n-1}]|_{X_{n-1}= x}= Q(x,dx'),
\end{equation}
for all $n\ge 1$, and with initial distribution $\mu_{0}$.

We also assume that given the state sequence $\{X_{n}, n\ge 0\}$, the observations $\{Y_{n}, n\ge 0\}$ are independent, the conditional distribution of $Y_{n}$ depends only on $X_{n}$ and that the conditional probability distribution $\mathbb{P}[Y_{n}\in dy | X_{n}= x]$ is absolutely continuous with respect to the Lebesgue  measure, \textit{i.e.}:
\begin{equation}\label{density}
 \mathbb{P}[Y_{n}\in dy | X_{n}= x]=\Psi(x,y)\,dy,
 \end{equation}
for some Borel measurable with respect to the couple \((x,y)\) function \(\Psi\).
The basic example which is to be covered will be the following:

We consider a discrete time filter for a Hidden Markov chain
$(X_n)$ with values in the Euclidean space $\mathbb{R}^{d}$, with
conditionally Markov observations $(Y_n)$ also from $\mathbb{R}^{1}$  satisfying the system
\begin{eqnarray}\label{e1ex}
&X_{n+1} = X_n+ b(X_{n}) + \xi_{n+1}, 
\quad (n\ge 0), \\\nonumber\\
&Y_n = h(X_n) + V_n \quad (n\ge 1),
\label{e2ex}
\end{eqnarray}

\noindent
where $(\xi_n,V_n)$ is a sequence of IID random vectors
of dimension $d+\ell$ with densities  $q_{\xi} (x),\, q_{V}(y)$,
$b(\cdot)$ is a $d$-dimensional vector-function,
$h(\cdot)$ an $\ell$-dimensional vector-function,
 that is,
\begin{equation}\label{Q}
Q(x,dx') =
q_{\xi}((x'-x-b(x)))\,dx',
\end{equation}
\begin{equation}\label{Psi}
\Psi(x,y)=q_{v}(y-h(x))
\end{equation}
(remind that $q_\xi$ and $q_{v}$ denote the densities  of $\xi_1$ and $V_{1}$ respectively).

The problem addressed in this paper is as follows. Assume that the exact parameters of the model \eqref{kernel}--\eqref{density} --- i.e.,  the initial distribution
$\mu_0$, the transition kernel $Q(x,dx')$ and the conditional density of the observations  $\Psi(x,y)$ --- are known with some errors or that we know only an approximations to the exact characteristics of this model. Hence, the statistician is unable to use the exact optimal filtering algorithm for estimation of  $X_{n}$ at each time $n$,  and he is left to apply a filtering algorithm with wrong parameters and with additional errors in the algorithm itself.

Under such conditions, the goal is to investigate the asymptotic behaviour  of this error in the
available algorithm in the long run. It follows from the earlier results on the subject -- see \cite{KV:08} --  that it is sufficient to work with errors  in the kernels assuming  that initial distribution $\mu_{0}$ is known exactly. (If not, it may be tackled by using the methods and results from \cite{KV:08}.)
More precisely the
setting will be explained in the section \ref{Main question} below.

Throughout the paper, we denote the wrong transition kernel and conditional density of the observations by $P(x,dx')$ and by $\Xi(x,y)$ respectively. Assumptions will be stated in the form of $Q, P$ and $\Psi, \Xi$ and examples will be provided in terms of the coefficients and properties of the original system (\ref{e1ex})--(\ref{e2ex}).

\subsection{Main result}\label{Main question}
To explain the main problem  addressed in this paper in detail, we should formulate exact and wrong filtering algorithms. Let us remaind to readers the exact filtering algorithm. The problem of nonlinear filtering is to compute at each time $n$, the conditional probability distribution $\mu_{n}$ of the state $X_{n}$ given  the observation sequence
$Y_{1:n}=\{Y_{1}, \, Y_{2}, \,\ldots \,Y_{n}, \,\}$, \textit{i.e.}:
$$
\mu_{n}(A)=\mu_{n}^Y (A)=\mathbb{P}[X_{n}\in A \mid Y_1, \ldots Y_n].
$$

Using Bayes' formula, the exact  posterior filtering conditional measure
can be represented as a probability measure for any $Y$ via the following random non-linear
operator $\bar S_n^{Y,\mu_0}$, applied to the initial measure $\mu_0$,
\begin{eqnarray}
\nonumber \mu_{n}(dx_n)=\mu_{n}^Y(dx_n)= \mathbb{P}_{\mu_0}(X_n \in dx_n \mid Y_1, \ldots Y_n)
 \\ \nonumber \\ \nonumber
= \int_{\mathbb{R}^{n}}
\prod_{i=1}^{n} Q(x_{i-1},dx_{i}) d_{i}^{\mu_0}\Psi(x_{i}, Y_i)
\mu_0(dx_0)
\\ \nonumber \\ \label{QR}
= \frac{1}{c_{n}^{\mu_0}} \int_{\mathbb{R}^{n}} \prod_{i=1}^{n} Q(x_{i-1},dx_{i}) \Psi(x_i,
Y_i) \mu_0(dx_0) =: \mu_0 \bar S_n^{Y,\mu_0} (dx_n).
\end{eqnarray}
Here 
$\Psi(x_i,y_i)$ is a conditional density of $Y_i$ at $y_i$, given
$X_i=x_i$, and $Q(x,dx')$ is a transition kernel for the Markov
chain $X_n,\, n\ge 0.$
The random
normalization constant $c_n^{\mu_0}$ is defined as follows,
$$
c_i^{\mu_0} =  \left.E_{\mu_0}\prod_{j=1}^{i} \Psi
(X_j,y_j) \right| _{y_1=Y_1, \ldots y_i=Y_i},
$$
and, correspondingly,
$$
d_i^{\mu_0} = \left. \frac{c_{i-1}^{\mu_0}}{c_{i}^{\mu_0}} =
\frac{E_{\mu_0}\left(\prod_{j=1}^{i-1} \Psi (X_j,y_j)\right) }
{E_{\mu_0}\left(\prod_{j=1}^{i} \Psi (X_j,y_j) \right) }
\right|_{y_1=Y_1, \ldots y_i=Y_i} .
$$
Hence,
$$
c_n^{\mu_0} = c_n^{Y,\mu_0} = \int_{\mathbb{R}^{n+1}} \prod_{i=1}^{n} Q(x_{i-1},dx_{i}) \Psi(x_i,
Y_i) \mu_0(dx_0).
$$

~

Now, the ``wrong filtering algorithm '' can be formulated more precisely as follows. Recall that it is assumed that we do not know
 the transition kernel $Q(x,dx')$ and the conditional density of the observations  $\Psi(x,y)$ exactly, but only some  approximations $P(x,dx')$ and  $\Xi (x,y)$ respectively.
Hence we can define another sequence of measures $(\mu^{'}_{n}(A))_{n\ge 1}$ as follows:
\begin{eqnarray}
\mu^{'}_{n}(dx_n)=\int_{\mathbb{R}^{n}}
\prod_{i=1}^{n} P(x_{i-1},dx_{i}) \tilde  d_{i}^{\mu_0}\tilde  \Xi(x_{i}, Y_i)
\mu_0(dx_0)
\\ \nonumber \\ \label{QR'}
= \frac{1}{\tilde c_{n}^{\mu_0}} \int_{\mathbb{R}^{n}} \prod_{i=1}^{n} P(x_{i-1},dx_{i})  \Xi(x_i,
Y_i) \mu_0(dx_0) =: \mu_0 \tilde S_n^{Y,\mu_0} (dx_n),
\end{eqnarray}
where the ``wrong'' normalizing constant $\tilde c_{n}^{\mu_0}$ can be defined as follows:
$$
\tilde c_n^{\mu_0} = \int_{\mathbb{R}^{d+1}} \prod_{i=1}^{n} P^{}(x_{i-1},dx_{i}) \Xi^{}(x_i,
Y_i) \mu_0(dx_0).
$$
The problem is to estimate the asymptotic behaviour of the difference:
$$
E_{\mu_0} \|\mu_n' - \mu_n \|_{_{TV}}.
$$
We may not hope that this discrepancy goes to zero as \(n\to
\infty\), but just that under certain conditions it may remain small for all values of \(n\).

\subsection*{Assumptions}
\begin{enumerate}
\item[(A1)] --- bounded (small) local perturbations ---

We assume that
$$
\ln(\sup_{x,\,x^{\prime},\,z,\,y} \frac{Q(x,\,d x^{\prime})\Psi(z,y)}{P(x,\,dx^{\prime})\Xi(z,y)}) + \ln(
\sup_{x,\,x^{\prime},\,z,\,y} \frac{P(x,\,d x^{\prime})\Xi(z,y)}{Q(x,\,dx^{\prime})\Psi(z,y)})=q < \infty.
$$

\item[(A2)] --- local mixing ---

We assume
that for any $R>0$
\begin{eqnarray}\label{emixden}
  C_R =: \sup_{D_R}\,\, \left(\frac{Q(x_0,dx')}
{Q(v_0,dx')},  \frac{P(x_0,dx')}
{P(v_0,dx')}\right) <\infty, 
\end{eqnarray}
with $D_R:= \{(x_0, v_0,  x'):\, |x_0|,|v_0|,|x'| \le R, \}$
\item[(A3)] --- positiveness of conditional densities ---
$$
\Psi(x,y)>0, \quad \Xi^{}(x,y)>0, \quad \forall \; x,y.
$$

\item[(A4)] --- condition of exponential recurrence in terms of the transition kernels:

there exist $q\in (0,1)$, $R, K, c >0$ such that for $|x|\ge R$,

\begin{equation}\label{eqrho}
\left(\int \exp(c|x'|) Q(x,dx')\right) \vee
\left(\int \exp(c|x'|) Q(x,dx')\right) \le \rho \exp(c|x|),
 \end{equation}
and for $|x|\ge R$,
\begin{equation}\label{eqK}
\left(\int \exp(c|x'|) Q(x,dx')\right) \vee
\left(\int \exp(c|x'|) Q(x,dx')\right) \le K.
 \end{equation}

~

This condition may be also re-written as follows:
$$
|a(x)|:= |\mathbb{E}^Q_x X_1| \le q|x| + K, \quad \& \quad |a'(x):=\mathbb{E}^{P}_x X_1| \le q|x| + K.
$$

\item[(A5)] --- uniformly small  influence of observations:
there exists $\delta>0$ such that
$$
\sup_{x,y,x'} \frac{\Psi(x,y)}{\Psi(x',y)}
\vee \sup_{x,y,x'} \frac{\Xi(x,y)}{\Xi(x',y)}
\le 1+\delta,
$$
and for $\rho$ from (A4)
$$
(1+\delta)\rho<1.
$$

Although the  assumption  (A4) requires explicitly exponential tails of the noise in the signal, we give examples with both exponential and polynomial tails showing in what situations the assumption (A5) may be verified. Polynomial examples may be useful in the future.

~

Sufficient conditions for the assumption (A4)  in terms of the equation (\ref{e1ex}) and its approximation
\begin{eqnarray}\label{e1ex2}
&X'_{n+1} = X'_n+ \tilde b(X'_{n}) + \xi'_{n+1},
\quad (n\ge 0)
, \\\nonumber\\
&Y'_n = \tilde h(X'_n) +  V'_n \quad (n\ge 1),
\label{e2ex2}
\end{eqnarray}
may be offered as follows.

\item[(A4')] --- another condition of exponential recurrence:
there exist $r > 0$ such that for $|x|$ large enough,
$$
|x + b(x)| \vee |x+ \tilde b(x)| \le |x| - r,
$$
and for any $R>0$
$$
\sup_{|x|>R}|x + b(x)| \vee \sup_{|x|>R}|x+ \tilde b(x)| < \infty,
$$
and
$$
\mathbb{E}\xi_1 = \mathbb{E}\xi'_1 = 0,
$$
and finally, there exists $\epsilon>0$ such that
$$
\mathbb{E}\exp(\epsilon |\xi_1|)
+ \mathbb{E}\exp(\epsilon |\xi'_1|)< \infty.
$$

\end{enumerate}

~

\begin{Remark}\label{rem0}
Apparently, (A4a) implies the following (notations are taken from the calculus in the proof -- see below): there exists a constant $C<\infty$ such that for any probability measure $\mu$,
$$
\sup_y \int\int e^{\epsilon |x'|} \mu  \bar Q^{y,\mu} (x,dx') \mu(dx) \le C  \int e^{\epsilon |x|} \mu (dx).
$$
Equivalently,  we could say that for any probability measure $\mu$ we have,
\begin{eqnarray}\label{a5}
\displaystyle \sup_\omega \mathbb{E}_\mu \left(e^{\epsilon |X_1|} | Y_1\right) \le C  \int e^{\epsilon |x|} \mu (dx).
\end{eqnarray}
The same holds true for the approximation  kernel $P$.
\end{Remark}

~

The following is the main result of the paper.
\begin{Theorem}\label{Thm1}
Let $\displaystyle \int e^{\epsilon |x|} \,\mu_0(dx) < \infty$. Then, under the assumption (A1) -- (A5) above, there exists a constant $C>0$ such that the following bound
holds true:

\begin{equation}\label{ethm2}
\sup_n \mathbb{E}_{\mu_0} \| \mu'_n-\mu_n \|_{_{TV}} \le Cq.
\end{equation}
\end{Theorem}

\begin{Remark}\label{remcon}
Assumption  (A1) is valid for example for the model \eqref{e1ex}-\eqref{e2ex} with
$q_{v}=q_{\xi}= C\exp (-|x|)$ and with $\tilde  b^{}: \tilde b^{}(x)=b(x)$ if $|x|>K$ for some $K$.

Note that the value of $q$ may be an arbitrary value greater than zero; however, the question is meaningful if this constant is small.
It would be also nice to localize this condition, see the Remark \ref{dream} in the sequel; we leave it till further studies.

\end{Remark}

\begin{Remark}
Let us show how the assumption (A5) may be checked. \\

{\it Example 1.} Recall that
$\Psi(x,y)=q_{v}(y-h(x))$. Let $q_v(y) = c (1+|y|)^{-m}, \, |y|\ge M$. Assume that $\sup_x |h(x)| \le \delta'$ and this $\delta'$ is small. Then for $|y|$ large enough,
\begin{eqnarray*}
\frac{\Psi(x,y)}{\Psi(x',y)} = \frac{q_{v}(y-h(x))}{q_{v}(y-h(x'))}
 \\\\
= \frac{(1+|y-h(x)|)^{-m}}{(1+|y-h(x')|)^{-m}} = \frac{(1+|y-h(x')|)^{m}}{(1+|y-h(x)|)^{m}}
 \\\\
= \left(\frac{1+|y-h(x')|}{1+|y-h(x)|}\right)^m
\le \left(\frac{1+|y|+\delta'}{1+|y|-\delta'}\right)^m
 \\\\
=  \left(\frac{1+\delta'/(1+|y|)}{1-\delta'/(1+|y|)}\right)^m
\approx 1 + 2m\delta'/(1+|y|)
\le 1 + 2m\delta'.
\end{eqnarray*}
While for $|y|$ small enough, say, $|y|\le M'$,
\begin{eqnarray*}
\frac{\Psi(x,y)}{\Psi(x',y)} = \frac{q_{v}(y-h(x))}{q_{v}(y-h(x'))}
\le \frac{\sup_{|y''-y|\le \delta'}q_{v}(y'')}{\inf_{|y'-y|\le \delta'}q_{v}(y')}.
\end{eqnarray*}
If $q_v>0$ everywhere and is continuous, then the latter right hand side is close to 1 uniformly in $|y|\le M'$.

~

{\it Example 2.} Let $q_v(y) \sim c (1+|y|)^{-m}, \, |y|\to \infty$. Assume that $\sup_x |h(x)| \le \delta'$ and this $\delta'$ is small. Then for $|y|\to\infty$,
\begin{eqnarray*}
\frac{\Psi(x,y)}{\Psi(x',y)} = \frac{q_{v}(y-h(x))}{q_{v}(y-h(x'))}
 \\\\
\sim \frac{(1+|y-h(x)|)^{-m}}{(1+|y-h(x')|)^{-m}} = \frac{(1+|y-h(x')|)^{m}}{(1+|y-h(x)|)^{m}}
 \\\\
= \left(\frac{1+|y-h(x')|}{1+|y-h(x)|}\right)^m
\le \left(\frac{1+|y|+\delta'}{1+|y|-\delta'}\right)^m
 \\\\
=  \left(\frac{1+\delta'/(1+|y|)}{1-\delta'/(1+|y|)}\right)^m
\approx 1 + 2m\delta'/(1+|y|)
\le 1 + 2m\delta'.
\end{eqnarray*}
For $y$ bounded, the ratio remains close to one if $\delta'$ is small enough.

~

{\it Example 3.} Let $g_v(y) = c\exp(-|y|)$ for $|y|\ge M$, and $g_v$ be  continuous. Assume that $h$ is bounded and small: $|h(y)| \le \delta'$. Then, for any value of $y$,
\begin{eqnarray*}
\frac{\Psi(x,y)}{\Psi(x',y)} = \frac{q_{v}(y-h(x))}{q_{v}(y-h(x'))}
 \\\\
= \frac{c\exp(-|y-h(x)|)}{c\exp(-|y-h(x')|)} = \exp(-|y-h(x)| + |y-h(x')|)
 \\\\
\le \exp(2\delta') \approx 1 + 2\delta'.
\end{eqnarray*}

~

{\it Example 4.} Let $g_v(y) \sim c\exp(-|y|)$ as $|y|\to\infty$ (equavalent, i.e., the ratio converges to one), and let $g_v$ be  continuous.
Then, for $|y|$ large enough we have,

\begin{eqnarray*}
\frac{\Psi(x,y)}{\Psi(x',y)} = \frac{q_{v}(y-h(x))}{q_{v}(y-h(x'))}
 \\\\
=  \frac{q_{v}(y-h(x))/c\exp(-|y-h(x)|)}{q_{v}(y-h(x'))/c\exp(-|y-h(x')|)}
  \frac{c\exp(-|y-h(x)|)}{c\exp(-|y-h(x')|)}
 \\\\
\sim \exp(-|y-h(x)| + |y-h(x')|)
\le \exp(2\delta') \approx 1 + 2\delta'.
\end{eqnarray*}
And for $|y|$ bounded the ratio remains close to one due to continuity.

~

{\it NB.} For Gaussian densities the assumption (A5) apparently fails. This could be possibly overcome if we manage to do some "localisation".

\end{Remark}

\begin{Remark}\label{dream}
It would be interesting to replace the Assumption  (A1) by a {\bf local} assumption of the type,
$$
\ln(\sup_{x,\,x^{\prime},\,z,\,y \in K} \frac{Q(x,\,d x^{\prime})\Psi(z,y)}{P(x,\,dx^{\prime}) \Xi(z,y)}) +
\ln(\sup_{x,\,x^{\prime},\,z,\,y \in K} \frac{P(x,\,d x^{\prime}) \Xi(z,y)}{Q(x,\,dx^{\prime}) \Psi(z,y)})=q_{K} < \infty
$$
with $q_K$ small, perhaps, in addition to
$$
\ln(\sup_{x,\,x^{\prime},\,z,\,y} \frac{Q(x,\,d x^{\prime})\Psi(z,y)}{P(x,\,dx^{\prime}) \Xi(z,y)}) +
\ln(\sup_{x,\,x^{\prime},\,z,\,y } \frac{P(x,\,d x^{\prime})\Xi(z,y)}{Q(x,\,dx^{\prime})\Psi(z,y)})=q < \infty,
$$
(with $q$ arbitrary and finite)
and to change the current statement of the  Theorem~\ref{Thm1} to the following one: the bound
holds true,
\begin{equation}\label{ethm22}
E_{\mu_0} \| \mu'_n-\mu_n \|_{_{TV}} < C q_{K}      + \ln\sup_{x\in K} \mathbb{E}_x\exp(\alpha\hat\tau),
\end{equation}
or, possibly,
\begin{equation}\label{ethm2a}
E_{\mu_0} \| \mu'_n-\mu_n \|_{_{TV}} < C q_{K}      +  q \ln\sup_{x\in K} \mathbb{E}_x\exp(\alpha\hat\tau),
\end{equation}
or likewise.
At the moment it is a conjecture that one of the bounds (\ref{ethm2}--\ref{ethm2a}) may hold true under less rigorous conditions than those in the Theorem \ref{Thm1}.

\end{Remark}

\subsection{Auxiliary results}
In  \cite{KV:08} it was proved, in particular,  that under the ``exponential'' assumptions equivalent to (A4) the following estimate holds true:
\begin{equation}\label{ethm222}
E_{\mu_0}\|\mu_0 \bar S_n^{Y,\mu_0} - \nu_0  \bar S_n^{Y,\nu_0} \|_{TV}
\le C \exp(-C'n).
\end{equation}

Here, we need some minor modification of  (\ref{ethm2}). Recall that the proof of this estimate was based on the inequalities (14) and (20) from \cite{KV:08}. In turn, (14)/\cite{KV:08} followed from (11) and (12)/\cite{KV:08}, while (20)/\cite{KV:08} was a corollary from the results about mixing for the recurrent and ergodic signal process. What is important for the present paper, is that the basic inequality (12) \cite{KV:08} admits an improved version under the condition that the initial Birkhoff distance (13) between measures \(\mu_0\) and \(\nu_0\) is finite:
\[
\rho(\mu_0, \nu_0) <\infty.
\]
In \cite{KV:08} this was not assumed and there was no reason for using such  an improved version; on the contrary, the absence of this assumption allowed to cover a wider class of processes.
However, now this will be important and the new version we need is as follows: 
\begin{equation}\label{from2}
\rho(\mu_n, \nu_n) \le \rho((\mu_n, \nu_n), (\nu_n,\mu_n)) \le C \pi_R^{k-1} \rho(\mu_0, \nu_0),
\end{equation}
with some $C>0$ and $\pi_R<1$.
We do not explain here what are exactly $k, \mu_n, \nu_n$, et al. because it would require to copy several pages from \cite{KV:08}, but use the notations from \cite{KV:08} verbatim.
The point is that as a result of these improvements, we now formulate a version of Theorem 1 from \cite{KV:08} as follows.

~

\noindent
{\bf Theorem 2 [version of \cite{KV:08}]} \\
{\em Let $\displaystyle \int e^{\epsilon |x|} \,\mu_0(dx) < \infty$. Then under the assumptions (A2) -- (A4), the following bounds
hold true: there exist constants $C, C(\mu), \alpha, \epsilon >0$ such that

\begin{equation}\label{ethm22a}
E_{\mu_0}\|\mu_0\bar S_n^{Y,\mu_0} - \nu_0\bar S_n^{Y,\nu_0} \|_{TV}
\le C(\mu_0) \exp(-\alpha n)\rho(\mu_0,\nu_0),
\end{equation}
with
$$
C(\mu) \le \int \exp(\epsilon |x|) \,\mu(dx).
$$
}

~

\noindent
Note that we do not assume (A1) here because the statement relates only to identical kernels: $P \equiv Q$, and also $\Xi  \equiv \Psi$, in which case (A1) holds automatically with $q=0$.

~

\noindent
Also note that both versions -- the Theorem 2 above and the Theorem \ref{Thm1} from \cite{KV:08} -- could be combined with the help of the value \(1\wedge \rho(\mu_0,\nu_0)\) in the right hand side.

~

\begin{Lemma}\label{Le0}
Under the assumption (A4'), the assumption (A4) holds true.
\end{Lemma}

\begin{Lemma}\label{Le1}
Assume that (A4) holds true.
Let  $\hat\tau \equiv\hat\tau_R=
\inf(t\ge 0:\, |X_t|\le R)$. Then
there exist \(C,c,K,\lambda>0\) such that
\begin{eqnarray}\label{e-mix}
 \mathbb{E}_x\exp(\lambda\hat\tau) \le
C\exp(c|x|),
\quad \sup_{t\ge 0}E_x e^{c |X_t|} \le
(K + e^{c |x|}),
 \nonumber \\ \\ \nonumber
 \mathbb{E}_x\exp(\lambda \hat\tau') \le C\exp(c|x|),
\quad
\sup_{t\ge 0}E_x e^{c |X'_t|} \le
(K + e^{c |x|}).
\end{eqnarray}
(Recall that the process \(X'\) is built via the kernel \(P\) and \(\hat\tau'\) via (\(X'\))).
\end{Lemma}

\begin{Lemma}\label{Le2}
Assume (A4) and (A5). Then
for the ``conditional kernels''
\[
\overline  Q^{y, \mu} (x,dx') :=
Q(x,dx') \frac{\Psi(x',y)}{\iint Q(x,dx_1) \Psi(x_1,y) \mu(dx)}\, ,
\]
and
\[
\overline  P^{y, \mu} (x,dx') :=
P(x,dx') \frac{\Xi(x',y)}{\iint P(x,dx_1) \Xi(x_1,y) \mu(dx)}\, ,
\]
the following bounds holds true: for $|x|$ large enough,
\begin{equation}\label{newq}
|\int e^{\epsilon |x'|} \bar  Q^{y, \mu}(x,dx')| \vee |\int e^{c |x'|} \bar  P^{y, \mu}(x,dx')| \le \rho' e^{c |x|},
\end{equation}
and  with any $R>0$ for $|x|\le R$,
\begin{equation}\label{newqR}
|\int e^{c |x'|} \bar  Q^{y, \mu}(x,dx')| \vee |\int e^{c |x'|} \bar  P^{y, \mu}(x,dx')| \le K',
\end{equation}
with $\rho' = \rho(1+\delta)$ and $K' = K(1+\delta)$. Moreover, since $\rho'<1$ (see the Assumption (A5)), the following inequality holds true,
\begin{eqnarray}\label{expy}
\fbox{$\displaystyle \sup_{t\ge 0} \mathbb{E}_x e^{c |X^{',Y}_t|} \le K' + e^{c |x|}.$}
\end{eqnarray}

\end{Lemma}

\section{Proofs of auxiliary Lemmata}


Proofs of Lemma \ref{Le0} and Lemma \ref{Le1} based on \cite{GV1} and \cite{GV2} with some changes will be added in the next version of this preprint.





~

{\em Proof of Lemma \ref{Le2}.}
We have due to (A5) and (A4),  if $|x|$ is large enough, then
\begin{eqnarray*}
E_{x} e^{c |X^{',Y}_{1}|}
= \int e^{c |x'|} \overline P^{y, \mu} (x,dx')
 \\\\
\le (1+\delta) \int e^{c |x'|} P^{} (x,dx')
\le (1+\delta) e^{c |x|}.
\end{eqnarray*}
If $|x|$ remains bounded, say, $|x|\le R$ with some $R>0$, then
\begin{eqnarray*}
E_{x} e^{c |X^{',Y}_{1}|}
= \int e^{c |x'|} \overline P^{y, \mu} (x,dx') \le (1+\delta) K = K'.
\end{eqnarray*}
Now we estimate,
\begin{eqnarray*}
E_{X^{',Y}_{n-1}} e^{c |X^{',Y}_{n}|} \le
1(|X^{',Y}_{n-1}| > R) \rho' e^{c |X^{',Y}_{n-1}|}
+ 1(|X^{',Y}_{n-1}|\le R) K'
 \\\\
= \rho' e^{c |X^{',Y}_{n-1}|} - 1(|X^{',Y}_{n-1}| \le R)  e^{c |X^{',Y}_{n-1}|}
+ 1(|X^{',Y}_{n-1}|\le R) K'
 \\\\
\le \rho' e^{c |X^{',Y}_{n-1}|} + K'. \hspace{3cm}
\end{eqnarray*}
Further,
\begin{eqnarray*}
E_{X^{',Y}_{n-2}} e^{c |X^{',Y}_{n}|} \le
E_{X^{',Y}_{n-2}}  \left(\rho' e^{c |X^{',Y}_{n-1}|} + K'\right)
 \\\\
\le K' + \rho' \left(\rho' e^{c |X^{',Y}_{n-2}|} + K'\right),
\end{eqnarray*}
and eventually by induction,
\begin{eqnarray*}
E_{x} e^{c |X^{',Y}_{n}|} \le
\frac{K'}{1-\rho'} + (\rho')^n  e^{c |x|} \le
\frac{K'}{1-\rho'} + e^{c |x|}.
\end{eqnarray*}
The Lemma \ref{Le2} is proved.

\section{Proof of Theorem \ref{Thm1}}

{\bf 1.} We will use the Birkhoff metric for positive measures,
see \cite{KLS:89}, and also \cite{AZ:97}, \cite{LGO:04} (where it is called
Hilbert metric; one more synonym is the projective metric),
\begin{equation}\label{Birkhoff}
\rho(\mu,\nu) = \left\{\begin{array}{ll} \ln \frac{\mbox{$(\inf s:
\; \mu \le s \nu)$}}{\mbox{$(\sup t: \; \mu \ge t\nu)$}},
& \mbox{if finite}, \\
+\infty, & \mbox{otherwise}.
\end{array}
\right.
\end{equation}
Another equivalent definition reads,
$$
\rho(\mu,\nu) = \left\{\begin{array}{ll} \ln \sup (d\mu/d\nu) +
\ln \sup (d\nu/d\mu), & \mbox{if finite}, \\
+\infty, & \mbox{otherwise}.
\end{array}
\right.
$$

For any measure $\mu$ we can define the following  nonlinear operator $\overline{S}_{k:n}$:
for $k<n$
$$
\mu \bS_{k:n} (A)=\mu \bS_{k:n}^{Y} (A)= c_{k:n} \int_{\bbR^{n-k}} \mathbf{1}(x_n\in A) \prod_{j=k+1}^n Q(x_{j-1},dx_j) \Psi(x_j,y_j) \mu(dx_k),
$$
with a normalizing constant $c_{k:n}$:
$$
c_{k:n}= c_{k:n}^{Y}=\left [ \int_{\bbR^{n-k}}     \prod_{j=k+1}^n Q(x_{j-1},dx_j) \Psi(x_j,y_j) \mu(dx_k)\right]^{-1}, \quad k<n,
$$
and for $k=n$ we let   $\mu \bS_{k:n}(A) =\mu(A)$.
Denote
$$
\mu \overline{Q}^y (A) := \frac{
\int_{\bbR\times \bbR} 1(x'\in A) \Psi(x,y) Q(x,dx') \mu(dx)
}{
\int_{\bbR\times \bbR} \Psi(x,y) Q(x,dx') \mu(dx)
}\,,
$$
and similarly \(\mu \overline{P}^{y}\) is defined,
\[
\mu \overline{P}^{y} (A) := \frac{
\int_{\bbR\times \bbR}  1(x'\in A) \Xi(x,y) P(x,dx') \mu(dx)
}{
\int_{\bbR\times \bbR} \Xi(x,y) P(x,dx') \mu(dx)
}\,.
\]
Now we have:
\begin{equation}\label{eq:1}
\mu'_n-\mu_n = \sum_{k=1}^n (\mu'_k \bS_{k:n} - \mu'_{k-1} \bS_{k-1:n}).
\end{equation}

\vspace{1cm}
It was noted in  \cite{LGO:04} that the following important property holds true:
\begin{equation}\label{ie}
\mu'_{k-1} \bS_{k-1:n} = (\mu'_{k-1} \overline{Q}) \bS_{k:n},
\end{equation}
where the operator \(\overline Q = \overline  Q^{Y_{k}}\) is defined as
\[
\mu \overline  Q^{Y_{k}} (dx') =
\frac{\int Q(x,dx') \Psi(x',Y_{k}) \mu(dx)}{\iint Q(x,dx_1) \Psi(x_1,Y_{k}) \mu(dx)}\, ,
\]
or,
\[
\overline  Q^{Y_{k}} (dx') =
Q(x,dx') \frac{\Psi(x',Y_{k}) }{\iint Q(x,dx_1) \Psi(x_1,Y_{k}) \mu(dx)}\, ,
\]

~

For the reader's convenience we recall the reasoning.  Indeed,
$$
\mu \bS_{k:n} = \frac{1}{c_{k:n}^{\mu}} \mu S_{k:n},
$$
where the linear non-normalized operator $S_{k:n}$ is defined as follows:
$$
\mu S_{k:n}(A) = \mu S_{k:n}^Y (A) = \int_{\bbR^{n-k}} \mathbf{1} (x_n\in A) \prod_{j=k+1}^{n}\Psi(x_j, Y_j) Q(x_{j-1}, dx_j)
$$
Therefore, we have,
\begin{multline}
\mu S_{k:n}(A) = \int_{\bbR^{n-k-1}} \mathbf{1}(x_n\in A) \prod_{j=k+2}^{n} \Psi(x_j,Y_j)Q(x_{j-1}, dx_j)  \cdot \\
\cdot \int_{\bbR} Q(x_k,dx_{k+1})  \Psi(x_{k+1},Y_{k+1}) \, \mu(dx_k) = \nu S_{k+1:n}(A),
\end{multline}
with the non-normalized measure $\nu(dx_{k+1})$ defined by the formula
$$
\nu(dx_{k+1}) = \int_{\bbR} Q(x_k,dx_{k+1})  \Psi(x_{k+1},Y_{k+1})  \mu(d x_k)
$$
Hence,
\begin{equation}\label{markpr}
\mu \bS_{k:n} = \frac{\mu S_{k:n}}{c_{k:n}^{\mu}} = \frac{\nu S_{k+1:n}}{c_{k:n}^{\mu}}
= \frac{\nu S_{k+1:n}}{c_{k+1:n}^{\nu}} \cdot \frac{c_{k+1:n}^{\nu}}{c_{k:n}^{\mu}}.
\end{equation}
Also note (follows from the calculus with $A=\bbR$) that
$$
c_{k:n}^{\mu}=c_{k+1:n}^{\nu},
$$
because
\begin{multline}
\mu S_{k:n}(\bbR^d) = \int_{\bbR^{d(n-k-1)}} \mathbf{1}(x_n\in \bbR^d) \prod_{j=k+2}^{n} \Psi(x_j,y_j)Q(x_{j-1}, dx_j)  \cdot \\
\cdot \int_{\bbR^d} Q(x_k,dx_{k+1})  \Psi(x_{k+1},y_{k+1}) \, \mu(dx_k) = \nu S_{k+1:n}(\bbR^d).
\end{multline}


The equation \eqref{markpr} implies that
$$
\mu \bS_{k:n} = \frac{\nu S_{k+1:n}}{c_{k+1:n}^{\nu}} = \frac{(\nu/\nu(\bbR)) S_{k+1:n}}{c_{k+1:n}^{\nu}/ \nu(\bbR)} = \frac{\widetilde{\nu}S_{k+1:n}}{c_{k+1:n}^{\widetilde{\nu}}},
$$
with
$$
\widetilde{\nu}(dx')= \widetilde{\nu}^{Y_k}(dx')=\frac{\nu(dx')}{\nu(\bbR)} = \mu \overline{Q}^{Y_k} (dx') = \frac{\int Q(x,dx') \Psi(x', Y_k) \mu(dx)}{\iint Q(x,dx_1) \Psi(x_1,Y_k) \mu(dx)}\, .
$$
So, indeed, the announced
important property (\ref{ie}) holds true.

~

Further, since $\mu'_n \bS_{n:n} = \mu'_n$ and $\mu_0 \bS_{0:n}=\mu_n$ and because $\mu'_0 =\mu_0$ and $\mu'_0 \bS_{0:n} = \mu_0 \bS_{0:n}~=~\mu_n$, we obtain,
\begin{equation}\label{eq:1a}
\mu'_n-\mu_n = \sum_{k=1}^n (\mu'_{k-1}\overline{P} \, \bS_{k:n} - (\mu'_{k-1} \overline{Q}) \bS_{k:n}),
\end{equation}
where \(\mu'_{k-1}\overline Q = \mu'_{k-1}\overline Q^{Y_{k}}\),
 \(\mu'_{k-1}\overline P = \mu'_{k-1}\overline P^{Y_{k}}\).
So, it follows that
\[
\|\mu'_n-\mu_n\|_{TV} \le
\sum_{k=1}^n \|\mu'_{k-1}\overline P \, \bS_{k:n} - (\mu'_{k-1} \overline{Q}) \bS_{k:n}\|_{TV},
\]
and
\begin{equation}\label{start}
E_{\mu_0} \|\mu'_n-\mu_n\|_{TV} \le
\sum_{k=1}^n  E_{\mu_0} \|\mu'_{k-1}\overline P^{Y_{k} } \bS_{k:n} - (\mu'_{k-1} \overline{Q}^{Y_{k} }) \bS_{k:n}\|_{TV}.
\end{equation}

~

\noindent
{\bf 2.} By virtue of the Theorem 2 
under our assumptions
we have,
\[
E_{\mu,\nu}\| \mu \bS_{0:n} - \nu \bS_{0:n} \|_{TV} \le C(\mu,\nu) e^{-\alpha n} \rho(\mu,\nu),
\]
where \(\alpha\) {\em does not depend on the initial measures}, while \(C(\mu,\nu)\) admits a bound

\[
C(\mu,\nu) \le
\int (e^{c |x|}\mu(dx) +  e^{c |x'|}\nu(dx')).
\]
with some \(c>0\).

~

Also recall that due to the Lemma \ref{Le1} (here the process \(X'\) corresponds to the kernel~\(P\)), \begin{eqnarray}\label{bd1}
\sup_{t\ge 0}
E_{\mu_0}e^{\epsilon |X_t|}
\le K + \int e^{\epsilon |x|} \, \mu_{0}(dx),
 \nonumber \\\\ \nonumber
\sup_{t\ge 0}
 E_{\nu_0}e^{\epsilon |X'_t|}
\le K + \int  e^{\epsilon |x'|}  \nu_0(dx'). \end{eqnarray}

~

\noindent
{\bf 3.}
Further, all of the above imply that
\begin{eqnarray}\label{ine2}
E_{\mu_k,\nu_k} \left(\|\mu_k \bS_{k:n}^{Y_{k+1}, \ldots, Y_n} - \nu_k \bS_{k:n}^{Y_{k+1}, \ldots, Y_n}\|_{TV} \mid Y_{1}, \ldots, Y_{k}\right)
 \nonumber \\ \nonumber \\ \nonumber
\le  C(\mu_k,\nu_k) e^{-\alpha (n-k)}\rho(\mu_k,\nu_k) \hspace{2cm}
 \\\\ \nonumber
\le  C e^{-\alpha (n-k)}\rho(\mu_k,\nu_k)
\,\int e^{c|x|}(\mu_k^{Y_{1}, \ldots, Y_k}(dx) + \nu_k^{Y_{1}, \ldots, Y_k}(dx)),
\end{eqnarray}
where the constants $C$ and $\alpha$ are non-random and do not depend on $k$. Denote
$$
D(\mu, \nu):=\int e^{c|x|}(\mu(dx) + \nu(dx)).
$$

By virtue of the inequalities \eqref{start} and (\ref{ine2}), we have,
\begin{eqnarray*}
E_{\mu_0} \| \mu'_n -\mu_n \|_{TV} \hspace{3cm}
 \\\\
\le \sum_{k=1}^{n} C \,E_{\mu_0} D(\mu'_k, \mu'_{k-1}\overline{Q}^{Y_{k})} e^{-\alpha (n-k)} \sup_\omega \rho(\mu'_k, \mu'_{k-1}\overline{Q}^{Y_{k}})
 \\\\
\le Cq \, \sum_{k=1}^{n} E_{\mu_0} D(\mu'_k, \mu'_{k-1}\overline{Q}^{Y_{k}}) e^{-\alpha (n-k)} .
\end{eqnarray*}
But $\mu'_k=\mu'_{k-1} \overline{P}^{Y_{k}}$ (the same operator but with the wrong kernel), so we obtain
$$
\rho(\mu'_k, \mu'_{k-1}\overline{Q}^{Y_{k}})\
=\rho(\mu'_{k-1}\overline{P}^{Y_{k}}, \mu'_{k-1}\overline{Q}^{Y_{k}}) \; \le q.
$$
Thus,
\begin{equation}\label{fin}
E_{\mu_0} \| \mu'_n -\mu_n \|_{TV} \le \sum_{k=1}^{n} Cq e^{-\alpha (n-k)} E_{\mu_0} D(\mu'_k, \mu'_{k-1}\overline{Q}^{Y_{k}}).
\end{equation}

~

\noindent
{\bf 4.}
It remains to estimate the term so as to show that it does not exceed some finite constant uniformly in \(k\). We have,
$$
E_{\mu_0} D(\mu'_k,\mu'_{k-1}\overline{Q}^{Y_{k}})
= E_{\mu_0} \int e^{c |x|} \mu'_k(dx) + E_{\mu_0} \int e^{c |x|} \mu'_{k-1}\overline{Q}^{Y_{k}}(dx).
$$
Here the first term in the right hand side satisfies the bound
$$
(\sup_k)\, E_{\mu_0} \int e^{c |x|} \mu'_k(dx)\le C \int e^{c |x|} \mu_0(dx) < \infty,
$$
according to the first part of (\ref{bd1}).

Let us inspect the second part of the right hand side. We have,
\begin{eqnarray*}
E_{\mu_0} \int e^{\epsilon |x|} \mu'_{k-1}\overline{Q}^{Y_{k}}(dx)
=
E_{\mu_0} \int e^{c |x_k|}\mu'_{k-1} \overline{Q}^{Y_{k}} (dx_k)
 \\\\
\stackrel{(A5)}{\le} (1+\delta)
E_{\mu_0} \int e^{c |x_k|}(\mu'_{k-1} Q) (dx_k)
 \\\\
\le K(1+\delta) + (1+\delta)E_{\mu_0} \int e^{c |x|}\mu'_{k-1}(dx)
 \\\\
\stackrel{(\ref{expy})}{\le} K' +  K' +(1+\delta) \int e^{c |x|}\mu_{0}(dx)<\infty,
\end{eqnarray*}
according to
(\ref{expy}) and the assumption (A5): see the Lemma \ref{Le2}. The Theorem \ref{Thm1} is proved.

\section*{Acknowledgements}
For the second author, the article was prepared within the framework of a subsidy granted to the HSE by the Government of the Russian Federation for the implementation of the Global Competitiveness Program. The same author is also grateful for support of the RFBR grant 13-01-12447-ofi-m2.

~

\end{document}